\documentclass{ifacconf}
\usepackage{natbib}  

\usepackage{graphicx}
\usepackage{caption}
\usepackage{tabularx,booktabs}
\usepackage{comment}
\usepackage{float}
\usepackage{color}
\usepackage{multicol}
\usepackage{amsmath,amssymb}
\usepackage{amsfonts}
\usepackage{textcomp}
\usepackage{psfrag}
\usepackage{multimedia}
\usepackage{fancybox}
\usepackage{mathtools}

\definecolor{wheat}{rgb}{0.96,0.87,0.70}

\usepackage{amsmath}

\DeclareMathOperator*{\argmin}{arg\,min}
\allowdisplaybreaks

\newcommand{\vect}[1]{\ensuremath{\boldsymbol{\mathrm{#1}}}}

\newtheorem{Lemma}{Lemma}
\newtheorem{Definition}{Definition}

\definecolor{dirk}{rgb}{0.3,1,0.3}

\definecolor{seb}{rgb}{1,1,0.5}

\renewcommand\refname{References and Notes}

\newcommand {\matr}[2]{\left[\begin{array}{#1}#2\end{array}\right]}

\newcommand {\cmatr}[2]{\left\{\begin{array}{#1}#2\end{array}\right.}

\newcounter{lastnote}

\begin{document} 
\begin{frontmatter}



\title{Data-Driven Predictive Control and MPC: \\ Do we achieve optimality?}
\author[First]{A.S. Anand}
\author[First]{S. Sawant}
\author[First]{D. Reinhardt}
\author[First]{S. Gros} 

\address[First]{Dept. of Engineering Cybernetics, Norwegian University of Science and Technology (NTNU), Norway. E-mail: sebastien.gros@ntnu.no}

\begin{abstract}
In this paper, we explore the interplay between Predictive Control and closed-loop optimality, spanning from Model Predictive Control to Data-Driven Predictive Control. Predictive Control in general relies on some form of prediction scheme on the real system trajectories. However, these predictions may not accurately capture the real system dynamics, for e.g., due to stochasticity, resulting in sub-optimal control policies. This lack of optimality is a critical issue in case of problems with economic objectives. We address this by providing sufficient conditions on the underlying prediction scheme such that a Predictive Controller can achieve closed-loop optimality. However, these conditions do not readily extend to Data-Driven Predictive Control. In this context of closed-loop optimality, we conclude that the factor distinguishing the approaches within Data-Driven Predictive Control is if they can be cast as a sequential decision-making process or not, rather than the dichotomy of model-based vs. model-free. Furthermore, we show that the conventional approach of improving the prediction accuracy from data may not guarantee optimality.
 
\end{abstract}

\end{frontmatter}

\section{Introduction}
Model Predictive Control (MPC) is a popular and mature method for the control of complex systems. A core principle in MPC is that the predictions supporting an MPC scheme are based on simulations using an internal model of the real system \citep{MPCbook}. This model can be based on a state-space description of the system or purely on input-output sequences. In addition, it can be based on first principles (white box), partially data-driven (gray box), or purely data-driven (black box) \citep{Mesbah2022}. The accuracy of the MPC predictions can be crucial if the closed-loop performance of the policy produced by the MPC scheme is important \citep{darby2012mpc}. This can be a challenge when the real system is stochastic. 

The introduction of Machine Learning (ML) into MPC has been intensively studied in the recent literature \citep{Maiworm2021, Mesbah2022}. The hope is arguably in part to alleviate the issues related to inaccurate models \citep{hewing2020learning}. A classic approach is to use ML techniques to improve the MPC model using real-world data, typically via ridge regression \citep{wu2019machine}. A core issue identified in these approaches is that they do not tie model improvement to the closed-loop performance of the MPC scheme. Rather, they construct a model that best fit the data, hoping that this will translate into improved closed-loop performance. However, this assumption is not well supported in theory and not always verified in practice \citep{gros2019data, piga2019performance}. 

In this paper, we propose sufficient conditions for the MPC model such that the MPC policy achieves optimal closed-loop performance. We then propose an analysis of these conditions, which provides strong insights into when a model that best fits the data in the expected value sense can lead to optimal closed-loop performance. We observe that these proposed conditions are intimately linked to the fact that MPC operates via sequential decision-making. As a result, they apply to all forms of MPC as long as the predictions are based on internal model simulations.

In contrast to MPC, Data-Driven Predictive Control (DDPC) does not build its predictions by simulating a model \citep{krishnan2021direct}. Instead, it builds predictions directly from past input-output data collected from the real system \citep{yang2015data, coulson2019data, dirk2023}. Bypassing the construction of a model to support the predictions is sometimes presented as beneficial on the basis that such direct data-driven predictions are then tailored to the optimization objectives \citep{coulson2019data}. However, this position is not well supported in theory nor in practice \citep{dirk2023, mattsson2023regularization, fiedler2021relationship}. 

More importantly, we observe that DDPC does not, in general, operate via sequential decision-making. In fact, the predictions built in DDPC do not necessarily achieve causality \citep{sader2023causality}. A crucial consequence of these observations is that the conditions for closed-loop optimality of MPC proposed in this paper do not readily extend to DDPC. This sets DDPC strikingly at odds with other Predictive Control approaches. However, we argue that in the context of closed-loop performance, the distinction between model-based vs. model-free Predictive Control is less relevant than the distinction between sequential decision-making or not. 

The rest of this paper is organized as follows: Sec. \ref{sec:background} introduces the necessary background knowledge,  Sec. \ref{sec:MPCasMDPmodel} presents the sufficient condition for closed-loop optimality of an MPC scheme. Sec. \ref{sec:DDPC} discusses its implication in the context of DDPC, and a discussion of the findings and conclusion are presented in Sec. \ref{sec:discussion}.

\section{Background}\label{sec:background}
\subsection{Markov Decision Processes}
In this paper, we will adopt the definition of closed-loop optimality in the context of discrete infinite-horizon discounted Markov Decision Processes (MDPs) \cite{sutton2018reinforcement}. The developments proposed in this paper can be readily extended to broader definitions, but we will remain within this framework for the sake of brevity. MDPs assume the system dynamics of the form,
\begin{align}
\label{eq:StateTransition}
\vect s_+ \sim \rho \left(\,.\,|\,\vect s, \vect a\,\right)
\end{align}
where $\vect s\in \mathbb S \subseteq \mathbb R^n $ and $\vect a \in \mathbb A \subseteq \mathbb R^m$ are a pair of state and action (or input), respectively, and $\vect s_+\in \mathbb S$ is the successor state. $\rho$ is a conditional probability density or, more generally, a conditional probability measure. For a given stage cost $\ell\,:\, \mathbb S\,\times\, \mathbb A \rightarrow \mathbb R$, the closed-loop performance of a policy $\vect \pi\,:\, \mathbb S\,\rightarrow\, \mathbb A$ is given by:
\begin{align}
\label{eq:MDPCost}
J\left(\vect\pi\right) = \mathbb E\left[\left.\sum_{k=0}^\infty \gamma^k \ell\left(\vect s_k,\vect a_k\right)\,\right |\, \vect a_k =\vect\pi\left(\vect s_k\right) \right],
\end{align}
where $\gamma\in (0,1]$ is a discount factor modeling the probabilistic lifetime of the system and ensuring the boundedness of the infinite-horizon problem. The expectation $\mathbb E$ is over the Markov Chain resulting from \eqref{eq:StateTransition} in closed-loop with policy $\vect\pi$. Solving the MDP consists in determining an optimal policy $\vect\pi^\star$ by minimizing \eqref{eq:MDPCost}, i.e. 
\begin{align}
\label{eq:OptPolicy}
\vect\pi^\star = \argmin_{\vect\pi\in\Pi}\, J\left(\vect\pi\right).
\end{align}
In this paper, we assume that \eqref{eq:StateTransition}, $\ell$ and $\Pi$ are such that \eqref{eq:MDPCost} is integrable. Under these assumptions, the solution to the MDP is characterized by the Bellman equations, 
\begin{subequations}
\label{eq:Bellman}
\begin{align}
Q^\star\left(\vect s,\vect a\right) &= \ell\left(\vect s,\vect a\right) + \gamma \mathbb E\left[ V^\star\left(\vect s_+\right)\,|\, \vect s,\vect a \right] \, ,\\
V^\star\left(\vect s\right) &= \min_{\vect a}\,\, Q^\star\left(\vect s,\vect a\right) \label{eq:Vstar:Bell} \, , \\
\vect\pi^\star\left(\vect s\right) &= \mathrm{arg} \min_{\vect a}\,\, Q^\star\left(\vect s,\vect a\right) \, . \label{eq:Pistar:Bell}
\end{align}
\end{subequations}
where $V\, \text{and } Q$ represents value and action-value functions respectively, and $V^{\star}\, \text{and } Q^{\star}$ are their optimal values. In the MDP framework, constraints on the state and action can be enforced by assigning infinite values to the stage cost $\ell$ for constraint violation. More specifically, if a stage cost $L\left(\vect s,\vect a\right)$ is to be minimized by the optimal policy $\vect\pi^\star$ while respecting a set of constraints $\vect h\left(\vect s,\vect a\right)\leq 0$, then
\begin{align}
 \ell\left(\vect s,\vect a\right) = L\left(\vect s,\vect a\right) + \cmatr{ccc}{0&\mathrm{if}&\vect h\left(\vect s,\vect a\right) \leq 0\\
\infty&\mathrm{if}&\vect h\left(\vect s,\vect a\right) > 0 } 
\end{align}
is used. We finally observe that equations \eqref{eq:Bellman} are generally very difficult to solve due to the curse of dimensionality of Dynamic Programming.

\subsection{Model Predictive Control}
\label{sec:EMPC}
An MPC scheme \citep{MPCbook} solves the following optimization problem, for a given state $\vect s_k$ at time $k$:
\begin{subequations}
\label{eq:MPC0}
\begin{align}
V^\mathrm{MPC}\left(\vect s_k\right)=\min_{\vect x,\vect u}&\quad \gamma^{N} T\left(\vect x_{N}\right) + \sum_{i=0}^{N-1} \gamma^k L\left(\vect x_i,\vect u_i\right)\label{eq:MPC0:Cost}\\
\mathrm{s.t.} &\quad \vect x_{i+1} = \vect f\left(\vect x_i,\vect u_i\right)\label{eq:MPC0:Dyn} \,,\\
&\quad \vect h\left(\vect x_i,\vect u_i\right)\leq 0 \label{eq:MPC0:Const} \,, \\
&\quad \vect x_{0} = \vect s_k,\quad \vect x_{N} \in \mathbb T \,. \label{eq:MPC0:Boundaries} 
\end{align}
\end{subequations}
Here $\vect f$ is the MPC prediction model, which can be in a classic state-space form, but it can also be an input-output prediction model if 
\begin{subequations}
 \label{eq:IO:State}
 \begin{align}
 \vect x_i &= \matr{c}{\vect y_{i}\ldots, \vect y_{i-p}, \vect u_{i-1},\ldots, \vect u_{i-l}}^\top,\\
  \vect s_k &= \matr{c}{\vect y_{k}\ldots, \vect y_{k-p}, \vect a_{k-1},\ldots, \vect a_{k-l}}^\top.
 \end{align}
 \end{subequations}
for some $p,l$ large enough to make $\vect s_k$ a Markov state, and where $\vect y$ are measured outputs of the system. The terminal set $\mathbb T$ and the terminal cost $T$ are typically used to account for the fact that problem \eqref{eq:MPC0} is solved over a finite horizon $N$, while the MDP is open-ended. We assume throughout this paper that,  $V^\mathrm{MPC} = \infty$ when \eqref{eq:MPC0} is infeasible. In this paper, we follow the MPC formulation with discounting \eqref{eq:MPC0}, but all the results presented in this paper are valid for the undiscounted case ($\gamma = 1$). 
Note that the choice $\gamma=1$ common in MPC places some specific requirements on the stage cost $\ell$ for \eqref{eq:MDPCost} to be well-posed. 

The solution of \eqref{eq:MPC0} delivers an input sequence $\vect u^\star_{0,\ldots,N-1}$ and the corresponding state sequence $\vect x^\star_{0,\ldots,N}$ that the system is predicted to follow, according to the MPC model $\vect f$. Because the model is typically inaccurate, this prediction is usually inexact. To address this, MPC \eqref{eq:MPC0} is solved at every discrete time $k$, using the latest state of the system $\vect s_k$. Importantly, only the first input of the sequence, i.e. $\vect u^\star_{0}$ is applied to the real system. Therefore, the MPC scheme  \eqref{eq:MPC0} produces a policy 
\begin{align}
\label{eq:MPC:Policy}
\vect\pi^\mathrm{MPC}\left(\vect s_k\right) = \vect u^\star_0
\end{align}
that assigns for every feasible state $\vect s_k$ a corresponding action $\vect u^\star_0$ implicitly defined by \eqref{eq:MPC0}. In this paper, we aim to clarify how an MPC solves the corresponding MDP \eqref{eq:MDPCost}-\eqref{eq:OptPolicy} and discuss sufficient conditions that the underlying predictive mechanism of the MPC scheme should satisfy to achieve $\vect\pi^\mathrm{MPC}=\vect\pi^\star$. 

\section{MPC as a solution for MDPs}\label{sec:MPCasMDPmodel}
The goal of this section is to propose formal conditions under which the MPC policy \eqref{eq:MPC:Policy} achieves closed-loop optimality, i.e. minimizes \eqref{eq:MDPCost} in the presence of stochasticity and model error. 
\subsection{MPC as an action-value function} 
In order to form conditions for the closed-loop optimality of MPC, we will detail how an MPC \eqref{eq:MPC:Policy} delivers an action-value function $Q^\mathrm{MPC}\left(\vect s,\vect a\right)$ and the conditions under which it corresponds to the optimal action-value function $Q^\star$ of the MDP. We propose the following definition of $Q^\mathrm{MPC}\left(\vect s,\vect a\right)$:
\begin{subequations}
\label{eq:MPC:Qmodel}
\begin{align}
Q^\mathrm{MPC}\left(\vect s_k,\vect a_k\right):=\min_{\vect x,\vect u}&\quad 
\eqref{eq:MPC0:Cost} \label{eq:MPC:Qmodel:Cost}\\
\mathrm{s.t.} &\quad \eqref{eq:MPC0:Dyn}-\eqref{eq:MPC0:Boundaries} \label{eq:MPC:Qmodel:const}\\ 
&\quad  \vect u_0=\vect a_k. \label{eq:MPC:Qmodel:Aconst}
\end{align}
\end{subequations}
We observe that \eqref{eq:MPC:Qmodel} is identical to MPC \eqref{eq:MPC0}, to the exception of constraint \eqref{eq:MPC:Qmodel:Aconst}, which assigns the initial control input of the MPC. This definition of $Q^\mathrm{MPC}$ is valid in the sense of fundamental Bellman relationships between optimal action-value functions, value functions, and policies, i.e. 
\begin{subequations}
\begin{align}
V^\mathrm{MPC}\left(\vect s\right) &= \phantom{\mathrm{arg}}\min_{\vect a}\, Q^\mathrm{MPC}\left(\vect s,\vect a\right),\\\quad \vect \pi^\mathrm{MPC}\left(\vect s\right) &= \mathrm{arg}\min_{\vect a}\, Q^\mathrm{MPC}\left(\vect s,\vect a\right)
\end{align}
\end{subequations}
hold by construction. This implies that MPC scheme \eqref{eq:MPC0} solves the MDP if the condition:
\begin{align}
\label{eq:MPC:PerfectModel}
Q^\mathrm{MPC}\left(\vect s,\vect a\right)=Q^\star\left(\vect s,\vect a\right),\quad \forall\,\vect s,\vect a
\end{align}
 holds, thereby delivering an optimal policy, the optimal value function, and the optimal action-value function of the MDP. In this case, MPC can provide a full representation of the MDP solution. Condition \eqref{eq:MPC:PerfectModel} can be relaxed such that the MPC policy matches the optimal policy $\vect\pi^\star$. Then $Q^\mathrm{MPC}$ needs to match $Q^{\star}$ only in the sense of
 \begin{align}
\label{eq:MPC:PerfectModel:argmin}
\mathrm{arg}\min_{\vect a}\,Q^\mathrm{MPC}\left(\vect s,\vect a\right)=\mathrm{arg}\min_{\vect a}\,Q^\star\left(\vect s,\vect a\right),\quad \forall\,\vect s.
\end{align}
Note that \eqref{eq:MPC:PerfectModel} implies \eqref{eq:MPC:PerfectModel:argmin}, but the converse is not true, hence making \eqref{eq:MPC:PerfectModel:argmin} less restrictive than \eqref{eq:MPC:PerfectModel}. While less restrictive,  \eqref{eq:MPC:PerfectModel:argmin} falls short of making the MPC scheme a fully satisfying solution of the MDP. Indeed, while an MPC scheme satisfying \eqref{eq:MPC:PerfectModel:argmin} would deliver the optimal policy, it would not predict the expected cost associated to a given state-action pair correctly. It is important to note that \eqref{eq:MPC:PerfectModel} can hold locally in the actions $\vect a$ for all states $\vect s$, and then deliver the optimal policy, the optimal value function and locally the correct action-value function. Hence, arguably, \eqref{eq:MPC:PerfectModel} holding locally in $\vect a$ is--for all practical purposes--as valuable as if it holds globally. 

An additional case of interest is if \eqref{eq:MPC:PerfectModel} is relaxed to be valid up to a constant, i.e. 
\begin{align}
\label{eq:MPC:PerfectModel:PlusConstant}
Q^\mathrm{MPC}\left(\vect s,\vect a\right) + Q_0 =Q^\star\left(\vect s,\vect a\right),\quad \forall\,\vect s,\vect a
\end{align}
for some constant $Q_0\in\mathbb R$. If the MPC scheme satisfies condition \eqref{eq:MPC:PerfectModel:PlusConstant}, then it delivers the optimal policy and the optimal value functions up to a constant. As we will elaborate in the next sections, this case is important to understand the connection between MPC and MDPs. Note that condition \eqref{eq:MPC:PerfectModel:PlusConstant} can be converted to condition \eqref{eq:MPC:PerfectModel} by simply adding a constant to the MPC costs \eqref{eq:MPC0:Cost} and \eqref{eq:MPC:Qmodel:Cost}.

\subsection{Conditions for MPC optimality} \label{sec:OptimalMPC0}
This section outlines the necessary conditions for an MPC scheme to fulfill \eqref{eq:MPC:PerfectModel:PlusConstant}, thereby delivering an optimal policy, along with a correct representation of the MDP solution. We provide the key condition in the following Theorem. \\
\begin{thm}
\label{Th:MPCOptimality}
If the MPC model $\vect f$ satisfies the equality:
\begin{align}
\label{eq:constantmodif}
\mathbb E\left[V^\star\left(\vect s_+\right)\,|\, \vect s,\vect a\, \right] - V^\star\left(\vect f\left(\vect s,\vect a\right)\right) = V_0
\end{align}
for a constant $V_0$, and $T = V^{\star}$ is used in MPC \eqref{eq:MPC:Qmodel} then \eqref{eq:MPC:PerfectModel:PlusConstant} holds.
\end{thm}
\vspace{.05cm}

\begin{pf} 
Consider the MPC \eqref{eq:MPC:Qmodel} and modified stage cost $\hat{L}(\vect s, \vect a) = L(\vect s, \vect a) + \gamma V_{0} $, and terminal cost, $T(\vect{s}) = V^{\star}(\vect{s})$
and the associated action value function
\begin{subequations}
\begin{align}
    \hat{Q}^{\mathrm{MPC}}(\vect s_k, \vect a_k) = &\min_{\vect x, \vect u} \, \gamma^N V^{\star}(\vect{x}_N) + \sum_{i=0}^{N-1} \gamma^i \hat{L}(\vect{x}_i, \vect{u}_i)) \label{eq:MPC:Qmodel2:cost}\\
    & \mathrm{s.t.}\quad \eqref{eq:MPC:Qmodel:const}, \eqref{eq:MPC:Qmodel:Aconst}
\end{align}    
\end{subequations}
Using \eqref{eq:constantmodif}, we get, 
\begin{align}
    \hat{Q}^{\mathrm{MPC}}(\vect s_k, \vect a_k) = Q^{\mathrm{MPC}}(\vect s_k, \vect a_k) + gV_0, \quad  g \in \mathbb R
\end{align}
where $g$ is some constant. 
\eqref{eq:MPC:Qmodel2:cost} can be further expanded as,
\begin{align}
    \eqref{eq:MPC:Qmodel2:cost} = & \min_{\vect x, \vect u} \, \gamma^N V^{\star}(\vect{x}_N) + 
    \sum_{i=0}^{N-1} \gamma^i \Big(L(\vect{x}_i, \vect u_i) \notag \\
    & + \gamma (\mathbb{E} [V^\star(\vect{x}_{i+1})| \vect{x}_i, \vect u_i)]
    - V^\star(\vect f(\vect{x}_i, \vect u_i)))\Big)\,.
\end{align}
Using a telescopic sum, we get, 
\begin{equation}\label{eq:Q and A}
\eqref{eq:MPC:Qmodel2:cost} = \min_{\vect x, \vect u} Q^{\star}(\vect s_k, \vect a_k) + \mathbb{E}\left[\sum_{i=1}^{N-1} \gamma^i A^{\star}\left(\vect{x}_i, \vect u_i\right)\right]
\end{equation}
wherein, $A^{\star}(\vect s, \vect a)$ denotes advantage function, defined as:
\begin{align}
A^{\star}(\vect s, \vect a)=\left\{\begin{array}{cc}
Q^{\star}(\vect s, \vect a)-V^{\star}(\vect s) & \text { if }\left|Q^{\star}(\vect s, \vect a)\right|<\infty \\
\infty & \text { otherwise }
\end{array}\right.
\end{align}
According to Bellman equalities, $\displaystyle\min_{\vect a} A^{\star}(\vect s, \vect a)=0$,
therefore, \eqref{eq:Q and A} reduces to $Q^\star$ and we get,
\begin{equation}
    Q^{\mathrm{MPC}}(\vect s, \vect a) + Q_0 = Q^{\star}(\vect s, \vect a), \quad  g \in \mathbb R
\end{equation}
where $Q_0 = g V_{0}$, which concludes the proof.
\end{pf}

The immediate observation on Theorem \ref{Th:MPCOptimality} is that it establishes a sufficient condition on the MPC model $\vect f$ to satisfy such that the MPC scheme is closed-loop optimal. This condition interlaces the optimal value function of the MDP, $V^\star$, and the prediction model $\vect f$ in a non-trivial way. Before further investigating the implications of Theorem  \ref{Th:MPCOptimality}, it is noteworthy that condition \eqref{eq:constantmodif} is \textit{unlike} the conventional criteria used in SYSID and ML for constructing a prediction model that best fits the data. This suggests that an MPC scheme using such conventional criteria to build $\vect f$ does not necessarily yield an optimal policy.  We will further detail that in the next section. 

In the context of MPC, regression-based approaches are generally trying to estimate a model $\vect f$  as an expectation over the system state transitions as follows,
\begin{align}
\label{eq:E:Fitting}
\vect f\left(\vect s,\vect a\right) = \mathbb E\left[\vect s_+\,|\,\vect s,\vect a \right],\quad \forall\, \vect s,\vect a \,.
\end{align}
Therefore, we will focus the rest of our discussions in the context of the family of models seeking to achieve \eqref{eq:E:Fitting}. In addition to their relevance within MPC, this choice is motivated by the fact that this family of models, in combination with Theorem \ref{Th:MPCOptimality}, produces interesting results regarding the closed-loop optimality of MPC. The subsequent sections will explore this in detail. 

Another interesting remark on the condition \eqref{eq:constantmodif} is that it applies only to the state transition model $\vect f$, which predicts in a one-step-ahead manner in MPC, and not to the entire prediction $\vect x_{1,\ldots,N}$ that is inherent in MPC decisions. This observation is closely related to the sequential decision-making nature of MPC, as attested in the proof of Theorem \ref{Th:MPCOptimality}. We will discus this remark later in the paper. 

\subsection{Local optimality of MPC with expected-value model} 
\label{sec:LocalOpt}
 In this section, we will investigate the consequence of Theorem \ref{Th:MPCOptimality} for models achieving \eqref{eq:E:Fitting}. If model $\vect f$ is constructed such that it achieves \eqref{eq:E:Fitting}, then condition \eqref{eq:constantmodif} becomes:
 \begin{align}
\label{eq:Delta:E}
\mathbb E\left[V^\star\left(\vect s_+\right)\,|\, \vect s,\vect a\, \right] - V^\star\left(\mathbb E\left[\vect s_+\,|\,\vect s,\vect a \right]\right) = V_0,\quad  \forall\, \vect s,\vect a,
\end{align}
i.e. it requires that the state transition expected value operator $\mathbb E\left[\cdot \right]$ commutes with the optimal value function $V^\star(\cdot)$ up to a constant $V_0$. We can stress here that condition \eqref{eq:Delta:E} is purely related to the properties of the MDP, i.e. the stochastic dynamics \eqref{eq:StateTransition}, and the MDP cost $\ell$ in \eqref{eq:MDPCost}, and not on the ones of the MPC scheme. 
In order to facilitate the discussion of \eqref{eq:Delta:E}, let us first consider the minimum attraction set $\mathbb F \subseteq \left\{\,\vect s\in\mathbb S\,|\, V^\star\left(\vect s\right)<\infty\right\}$ associated to dynamics \eqref{eq:StateTransition} and optimal policy $\vect \pi^\star$ implicitly defined as the minimum compact set such that
\begin{align}
\mathbb P\left[\,\vect s_+ \in \mathbb F\,|\, \vect s \in\mathbb F,\,\vect a\in\mathcal B(\vect\pi^\star(\vect s),r)\right] = 1
\end{align}
holds, where $\mathcal B(\vect\pi^\star(\vect s),r)$ is a ball of radius $r$ centred at $\vect\pi^\star(\vect s)$. Set $\mathbb F$ then defines the set where the MDP trajectories eventually converge, with the addition of a small set of action disturbances. We ought to stress that a compact set $\mathbb F$ does not necessarily exist. Indeed, set $\mathbb F$ is strongly related to the dissipativity of the MDP. For a dissipative MDP, the closed-loop state distribution under optimal policy $\vect\pi^\star$ converges in the sense of a dissimilarity measure, like KL-divergence, Wasserstein distance, total variation, to the optimal steady-state distribution $\bar\rho$ defined by \citep{gros2022economic}:
\begin{subequations}
\begin{align}
\min_{\bar \rho,\bar{\vect\pi}}&\quad \mathbb E_{\vect s\sim \bar \rho}\left[\ell\left(\vect s,\bar{\vect\pi}\left(\vect s\right)\right)\right]\\
\mathrm{s.t}&\quad \bar\rho\left(\cdot\right) = \int \bar \rho\left(\vect s_+\right)\rho\left(\vect s_+|\vect s,\bar{\vect\pi}\left(\vect s\right)\right)\mathrm d \vect s_+.
\end{align}
\end{subequations}
In particular, for $r=0$, set $\mathbb F$ is given by the support of $\bar\rho$ if that support is compact. Let us define:
\begin{align}
\label{eq:Delta} 
&\Delta\left(\vect s,\vect a\right):=\mathbb E\left[V^\star\left(\vect s_+\right)\,|\, \vect s,\vect a\, \right] - V^\star\left(\vect f\left(\vect s,\vect a\right)\right) 
\end{align}
with $\vect f\left(\vect s,\vect a\right)$ given by \eqref{eq:E:Fitting}. We observe that satisfying \eqref{eq:Delta:E} requires that $\Delta$ is constant. We then propose the following Lemma, which provides a useful analysis of $\Delta$. For simplicity, we will focus on the case where the dynamics \eqref{eq:StateTransition} has a compact support.
\vspace{.05cm}
\begin{Lemma} 
\label{Lem:LocalOptimality}
Let us assume that $V^\star$ is $N\geq 3$ times continuously differentiable on $\mathbb F$. 

Then for any $ \vect s\in\mathbb F$ and $\vect a\in\mathcal B(\vect\pi^\star(\vect s),r)$:
\begin{align}
\label{eq:Delta} 
\Delta\left(\vect s,\vect a\right)=&\frac{1}{2}\mathrm{Tr}\left(\Sigma\left(\vect s,\vect a\right)\nabla^2V^\star\left(\vect f\left(\vect s,\vect a\right)\right)\right) + R\left(\vect s,\vect a\right)
\end{align}
where $\Sigma\left(\vect s,\vect a\right)$ is the conditional covariance of the state transition \eqref{eq:StateTransition} and, using the multi-index notation, $R$ is bounded by:
\begin{align}
\label{eq:RBound}
\left|R\left(\vect s,\vect a\right)\right| \leq c\mu_N+ \sum_{|\vect\alpha|=3}^{N-1}\frac{1}{\vect\alpha!}D^{\vect\alpha} V^\star\left(\vect f\left(\vect s,\vect a\right)\right) \mu_{|\vect\alpha|}
\end{align}
for some constant $c\geq 0$, with the moments
\begin{align}
\mu_k\left(\vect s,\vect a\right) &= \mathbb E\left[\left\|\vect s_+-\vect f\left(\vect s,\vect a\right)\right\|^k\right].
\end{align}
\end{Lemma}
\begin{pf}
Under the proposed assumptions, the optimal value function $V^\star$ admits a $N^\mathrm{th}$-order Taylor expansion on $\mathbb F$ at $\vect f\left(\vect s,\vect a\right)$, i.e.
\begin{align}
&V^\star\left(\vect s_+\right) =  V^\star\left(\vect f\left(\vect s,\vect a\right)\right) + \left(\vect s_+-\vect f\left(\vect s,\vect a\right)\right)^\top   \nabla V^\star\\
&+\frac{1}{2} \left(\vect s_+-\vect f\left(\vect s,\vect a\right)\right)^\top \nabla^2 V^\star\left(\vect s_+-\vect f\left(\vect s,\vect a\right)\right) + R\left(\vect s,\vect a\right)\nonumber
\end{align}
holds for $\vect s_+\in\mathbb F$, where $ \nabla V^\star$ and $\nabla^2 V^\star$ are evaluated at $\vect f\left(\vect s,\vect a\right)$, and
$R\left(\vect s,\vect a\right)$ is the Taylor expansion of $V^\star$ at $\vect f\left(\vect s,\vect a\right)$ from order 3 to $N-1$, together with the remainder of order $N$. It follows that \eqref{eq:RBound} 

holds for constant $c$, given by:
\begin{align}
\label{eq:TaylorConstant}
c &= \frac{1}{N!}\max_{|\vect \alpha|= N,\,\vect s_+\in\mathbb F}\left|D^{\vect\alpha} V^\star\left(\vect s_+\right)\right|
\end{align}
By assumption, $\vect s_+\in\mathbb F$ if $\vect s\in\mathbb F$ and $\vect a\in\mathcal B(\vect\pi^\star(\vect s),r)$, which concludes the proof.
\end{pf}
Note that Lemma \ref{Lem:LocalOptimality} does not require a compact set $\mathbb F$ if $V^\star$ is smooth everywhere and equals its Taylor series. Then Lemma \ref{Lem:LocalOptimality} applies if the moments $\mu_k$ decay fast enough with $k$. We will not detail this case here. Now, we will discuss the implications of Lemma \ref{Lem:LocalOptimality} on the satisfaction of  \eqref{eq:constantmodif}. A fully formal discussion is not provided here for the sake of brevity and will be the object of future publications.
\vspace{-.2cm}
\subsubsection{Deterministic MDPs} are an obvious class readily satisfying \eqref{eq:constantmodif}. In that case $\Sigma = 0$ and $R=0$ such that $\Delta = 0$. While of limited interest, this observation shows that Theorem \ref{Th:MPCOptimality} and Lemma \ref{Lem:LocalOptimality} support the fact that a deterministic MDP can be solved by an MPC having the correct model, terminal cost and terminal set. 
\vspace{-.2cm}
\subsubsection{LQR problems} are a second case of interest. If the optimal value function $V^\star$ is purely quadratic on $\mathbb F$, and $\Sigma$ is constant, then we observe that $R=0$, and that $\nabla^2V^\star$ is constant. As a result, $\Delta\left(\vect s,\vect a\right)=\frac{1}{2}\mathrm{Tr}\left(\Sigma\nabla^2V^\star\right) = V_0$, is constant and \eqref{eq:constantmodif} is satisfied. This case corresponds to the class of LQR problems with an additive process noise independent of $\vect s$ and $\vect a$, and generally to the class of MDPs that are equivalent to an LQR on their attraction set $\mathbb F$. While these are arguably very specific cases, they inform us that an MDP that is locally equivalent to an LQR problem can be locally solved by an MPC scheme with an expectation-based model \eqref{eq:E:Fitting}. This observation is further discussed next.
\vspace{-.2cm}
\subsubsection{Local optimality of MPC:} Let us now discuss a more interesting consequence of Lemma \ref{Lem:LocalOptimality}. First we observe that if $\Sigma\left(\vect s,\vect a\right)$ and $\nabla^2V^\star\left(\vect f\left(\vect s,\vect a\right)\right)$ are bounded and Lipschitz continuous on $\mathbb F$, then the first term in \eqref{eq:Delta} is also Lipschitz continuous, and its Lipschitz constant is related to how strongly the covariance $\Sigma$ of the dynamics \eqref{eq:StateTransition} and the Hessian $\nabla^2V^\star$ change over $\mathbb F$, and on how large $\Sigma$ and $\nabla^2V^\star$ are over $\mathbb F$. We further observe that $R$ in \eqref{eq:Delta} is small or nearly constant if the moments $\mu_{k>2}$ are small or constant, and if the higher-order derivatives of $V^\star$ are small or nearly constant. Hence $\Delta$ in \eqref{eq:Delta} is close to constant on $\mathbb F$ if the MDP is sufficiently smooth and if
\begin{enumerate} 
\item $\mathbb F$ is small, i.e. if the MDP is dissipative, and converging to a narrow steady-state distribution $\bar\rho$, or
\item $\Sigma$, $\nabla^2V^\star$, and $R$ are nearly constant over $\mathbb F$, i.e. if the moments higher than 2 of the state transition $\rho$ in \eqref{eq:StateTransition} are small or if they do not depend strongly on $\vect s,\vect a$ and if the Hessian of the optimal value function does not vary much over $\mathbb F$.
\end{enumerate}

These observations are very technical but point to important practical observations. First, the class of tracking problems with smooth dynamics of compact support and where the stage cost $L$ is quadratic, can fall into the categories listed above. Indeed, tracking problems are typically dissipative, the attraction set $\mathbb F$ can be small, and the value function is smooth where it is bounded. For this class of problems, one can expect that an MPC model based on the expected value model \eqref{eq:E:Fitting} achieves near-optimal closed-loop performances at the steady-state distribution $\bar\rho$. 

Similarly, the class of tracking problems arguably extends to \textit{economic} MDPs with smooth dynamics and generic smooth stage costs $L$, but which achieve dissipativity. Similarly to tracking problems, dissipative MDPs can have a small attraction set $\mathbb F$ and can fall into the categories enumerated above. Therefore, for a smooth dissipative MDP, characterized by steady-state stochastic trajectories converging to a small set, an MPC scheme utilizing the expectation-based model in \eqref{eq:E:Fitting} delivers a policy that is locally closed-loop optimal.

For cases outside these classes of problems, an MPC scheme based on \eqref{eq:E:Fitting} can not guarantee optimality, even locally. Such cases can be, for e.g., MDPs not achieving dissipativity and MDPs based on non-smooth cost functions or non-smooth dynamics. In particular, MDPs subject to exogenous disturbances, such as the strong variation of energy prices in energy-related applications. 

\section{Data-Driven Predictive Control (DDPC)}\label{sec:DDPC}
In this section, we will explore the implications of the proposed results in the context of this in the DDPC. A DDPC scheme can be written as \citep{coulson2019data}:
\begin{subequations}\label{eq:DDPC}
\begin{align}
    V^\mathrm{PC}\left(\vect s_k\right) &= \min_{\vect y,\vect u} \quad \gamma^{N} T\left(\vect x_{N}\right) + \sum_{i=0}^{N-1} \gamma^i L\left(\vect y_{i},\vect u_{i}\right) \label{eq:DDPC:IO:Cost} \\
    \mathrm{s.t.} &\quad \vect \phi\left(H_{k},\matr{c}{\vect y_{k-p}\\\vdots\\ 
\vect y_{k+N}},\matr{c}{\vect u_{k-l}\\\vdots\\ \vect u_{k+N-1}}\right)=0 \label{eq:DDPC:Pred}\\
&\quad \vect h\left(\vect y_{i},\vect u_{i}\right)\leq 0 \label{eq:DDPC:IO:Const} \\
&\quad \vect x_{0} = \vect s_k,\quad \vect x_{N} \in \mathbb T \,.
\end{align}
\end{subequations}
Here, the stage cost $L$ and constraints $\vect h$ are typically defined in terms of the output $ \vect y$ alone. We use the notation $\vect x_0, \, \vect s_k$ as defined in \eqref{eq:IO:State}. Set $H_{k}$ collects past input-output trajectories from the system available at time $k$, that is,
\begin{align}
H_{k} = \left\{
\matr{c}{\vect y_{j-p}\\\vdots\\ \vect y_{j+N}},\,\matr{c}{ \vect u_{j-l}\\\vdots\\ \vect u_{j+N-1}}\quad \mathrm{s.t.}\quad j \leq k-N
\right\}
\end{align}
and $\vect \phi$ is a multi-step predictor that relates via condition \eqref{eq:DDPC:Pred} the past data $H_{k}$, the recent input-output data $\vect y_{k-p,\ldots, k}$, $\vect u_{k-l,\ldots, k-1}$, and the input-output predictions $\vect y_{k+1,\ldots, k+N}$ and $\vect u_{k,\ldots, k+N-1}$. Note that the former are provided to the DDPC scheme \eqref{eq:DDPC}, while the input-output predictions are decision variables labeled as $\vect y,\vect u$.

The multi-step predictor $\vect \phi$ takes different forms in the literature, and the consensus on which form can be deemed as \textit{model-based} vs. \textit{model-free} is arguably debatable.
One approach is to build $\vect \phi$ as a solution to problem:
\begin{subequations}
\label{eq:FL:Prediction}
\begin{align}
&\min_{\vect\alpha,  \vect y_{k+1,\ldots, k+N}}\quad \|\vect \alpha\|^2 \label{eq:FL:Regularization} \\
&\mathrm{s.t.}\quad \sum_{ j \leq k-N} \vect \alpha_j \matr{c}{\vect y_{j-p}\\\vdots\\ \vect y_{j+N} \\ \vect u_{j-l}\\\vdots\\ \vect u_{j+N-1}} = \matr{c}{\vect y_{k-p}\\\vdots\\ \vect y_{k+N} \\ \vect u_{k-l}\\\vdots\\ \vect u_{k+N-1}} \label{eq:FL:Constraints}
\end{align}
\end{subequations}
for any given future input sequence $\vect u_{k,\ldots,k+N-1}$, recent input-output data $\vect y_{k-p,\ldots,k}$, $\vect u_{k-l,\ldots,k-1}$, and past input-output data $H_k$. Problem \eqref{eq:FL:Prediction} then delivers an output prediction $ \vect y_{k+1,\ldots, k+N}$ corresponding to the given input sequence $\vect u_{k,\ldots,k+N-1}$. An important observation is that because \eqref{eq:FL:Prediction} is a linear least-squares problem, it takes an explicit prediction which delivers the output prediction $\vect y_{k+1,\ldots, k+N}$ through a linear multi-step predictor akin to a subspace system identification method, i.e.  \eqref{eq:FL:Prediction} takes an explicit solution
\begin{align}
\label{eq:SS:Pred}
\matr{c}{\vect y_{k+1}\dots \vect y_{k+N}}^\top   = \Psi  \matr{c}{\vect y_{k-p}\dots \vect y_{k} \,, \vect u_{k-l}\dots \vect u_{k+N-1}}^\top  
\end{align}
for a matrix $\Psi$ that can be defined via a simple Least-Squares regression on the past data $H_k$ \citep{mattsson2023regularization, shambhu2023, dirk2023}. 
Another approach (as in DeePC) consists in directly embedding \eqref{eq:FL:Prediction} in the Predictive Control problem \eqref{eq:DDPC}, see \citep{ coulson2019data}. In that context, the decision variables $\vect\phi$ become part of problem \eqref{eq:DDPC}, the regularization objective \eqref{eq:FL:Regularization} is weighted and added to the control objective \eqref{eq:DDPC:IO:Cost}, and constraint \eqref{eq:DDPC:Pred} is replaced by \eqref{eq:FL:Constraints}. In that context, a Lasso penalty is often used for \eqref{eq:FL:Regularization}. 


\subsection{DDPC as a solution to MDPs}
Let us now explore under what conditions a DDPC scheme \eqref{eq:DDPC} achieves closed-loop optimality. An action-value function $Q^\mathrm{PC}\left(\vect s_k,\vect a_k\right)$ can be built from \eqref{eq:DDPC} with the addition of constraint \eqref{eq:MPC:Qmodel:Aconst}. Therefore, similarly to an MPC scheme a DDPC scheme can be construed as a model of an MDP solution. Before proceeding, let us introduce a simple definition to facilitate the following discussion.
 \vspace{.02cm}
\begin{Definition} Let us label a predictive control scheme as \textit{self-consistent} if its value function $V^\mathrm{PC}$ satisfies a Bellman equation. In the case of \eqref{eq:DDPC}, self-consistency requires that there is a function $\vect f\left(\vect s,\vect a\right)$ such that the following holds,  
\begin{align}
V^\mathrm{PC}\left(\vect s\right) =\min_{\vect a}\,\, \ell\left(\vect s,\vect a\right) + \gamma V^\mathrm{PC}\left(\vect f\left(\vect s,\vect a\right)\right)\quad \forall \vect s\in \mathbb S\,. 
\end{align}
\end{Definition} 
  \vspace{.02cm}
Self-consistency essentially requires that a predictive control scheme can be cast as a time-invariant sequential decision-making process, in line with an infinite-horizon MDP \eqref{eq:MDPCost}-\eqref{eq:OptPolicy}. MPC schemes \eqref{eq:MPC0} is self-consistent by construction for an adequate choice of terminal cost $T$ and set $\mathbb T$. However, a DDPC in the form \eqref{eq:DDPC} is not necessarily self-consistent because its prediction scheme is not based on the simulation of an underlying dynamic system. Indeed, unlike for \eqref{eq:MPC0}, function $\vect f$ is not explicitly available for a DDPC scheme, and does not exist in general. Let us discuss the implication of this observation on the closed-loop optimality of DDPC. A DDPC scheme delivers an optimal policy if \eqref{eq:MPC:PerfectModel:PlusConstant} holds, i.e. if
\begin{align}
\label{eq:QMatchDDPC}
Q^\mathrm{PC}\left(\vect s,\vect a\right) = Q^\star\left(\vect s,\vect a\right) + Q_0,\quad \forall\vect s,\vect a
\end{align}
holds at least locally in $\vect a$, or if the minimizers in $\vect a$ of $Q^\mathrm{PC}$ and $Q^\star$ match for all $\vect s$, as in \eqref{eq:MPC:PerfectModel:argmin}. Unlike for \eqref{eq:MPC0} condition \eqref{eq:QMatchDDPC} for DDPC is not a statement on the its predictions alone, and there is no reason to expect that it is satisfied by DDPC schemes, even in the context in Sec. \ref{sec:LocalOpt}.

This leads to an interesting observation in the context of this paper, i.e. the key distinction between the state-based MPC schemes \eqref{eq:MPC0} and the DDPC scheme \eqref{eq:DDPC} is not \textit{model-based vs. model-free}, but rather that the latter does not form predictions by simulating an underlying dynamic system, and can therefore not necessarily be cast as a sequential decision-making process.

If \eqref{eq:DDPC} is self-consistent, then the discussions of Sec. \ref{sec:OptimalMPC0} apply in the context of DDPC. Self-consistency can then be promoted or even enforced when using a data-driven predictor in the explicit form \eqref{eq:SS:Pred}. This can be achieved by enforcing self-consistency via block-diagonal regularization or constraints in the regression forming $\Psi$, promoting a partial block-Toeplitz structure in matrix $\Psi$, see \citep{shambhu2023}. In contrast, self-consistency is difficult to promote in an implicit predictor \eqref{eq:FL:Prediction} or when embedding \eqref{eq:FL:Prediction} in \eqref{eq:DDPC}. Therefore, in its generic formulation, we can not argue about the optimality of DDPC as we did in the case of MPC. Our observation is that DDPC may have to adhere to a restrictive structure, as in  \eqref{eq:SS:Pred}, on its predictor to attain local optimality.

%
%

\section{Discussion \& Conclusions} \label{sec:discussion}
In this paper, we have discussed the closed-loop optimality of Predictive Control schemes, as viewed in the framework of Markov Decision Processes, i.e. when the system to be controlled is stochastic. The question holds particular relevance for \textit{economic} problems, where closed-loop performance is of primary importance, in contrast to \textit{tracking} problems, where it can be secondary.

Our investigation identified sufficient conditions that allow a Predictive Control scheme to attain closed-loop optimality. These conditions (Theorem \ref{Th:MPCOptimality}) challenge the conventional approach of improving the prediction accuracy to achieve the best closed-loop performance for a Predictive Control scheme. We then provide an analysis that reconciles this conventional approach with our conditions, but this is limited to smooth tracking (or dissipative) problems and only valid locally around the (stochastic) steady state of the problem. For non-smooth or non-dissipative problems, the best-fitting model does not necessarily allow the Predictive Control scheme to deliver an optimal policy, even locally.


A key observation from our investigation across different Predictive Control schemes is that the proposed sufficient conditions are only applicable if the Predictive Control scheme can be viewed as a sequential decision-making process. This requires that the predictions are based on an underlying one-step ahead simulation, regardless of whether it is performed explicitly or not. We have labeled this property \textit{self-consistency}.

While MPC is self-consistent by construction,  the same does not hold true for DDPC in its direct form. Hence, the optimality conditions we have identified are not automatically applicable to DDPC. Therefore, to answer whether DDPC achieves optimality, we do not know unless it adheres to the restrictive structure proposed in \eqref{eq:SS:Pred}. However, our observations suggest that--in the context of closed-loop performance--the most important distinction between different Predictive Control methods is not between model-based and model-free approaches but rather between approaches that are self-consistent or not.

Finally, we note that the observations on local closed-loop optimality of MPC pertain to models seeking expected-value predictions \eqref{eq:E:Fitting}. While expected-value prediction is the objective of several SYSID and ML approaches, other forms of predictions, such as Maximum Likelihood Estimation, Quantile Regression, and Stochastic Predictions, are also used. However, we do not expect these approaches to have strong properties with respect to Theorem \ref{Th:MPCOptimality}. However, we leave it for further work, a deeper investigation of various prediction techniques in the context of Theorem \ref{Th:MPCOptimality}.

\bibliography{bib} 

\end{document}